\documentclass[10pt,reqno]{amsart}
\usepackage{amsthm,amsfonts,amssymb,euscript}

\def\normo#1{\left\|#1\right\|}

\def\norm#1{\|#1\|}
\def\jb#1{\langle#1\rangle}

\def\wh#1{\widehat{#1}}

\newcommand{\N}{{\mathbb N}}

\newcommand{\R}{{\mathbb R}}
\newcommand{\C}{{\mathbb C}}
\newcommand{\Z}{{\mathbb Z}}
\newcommand{\ft}{{\mathcal{F}}}

\newcommand{\les}{{\lesssim}}
\newcommand{\ges}{{\gtrsim}}

\newcommand{\Sch}{{\mathcal{S}}}
\newcommand{\supp}{{\mbox{supp}}}

\theoremstyle{plain}
  \newtheorem{theorem}[subsection]{Theorem}
  
  \newtheorem{lemma}[subsection]{Lemma}

\theoremstyle{remark}

\theoremstyle{definition}

\numberwithin{equation}{section}

\begin{document}
\title[Quadratic Schr\"odinger equation]{Remark on Well-posedness of
  Quadratic Schr\"odinger equation with nonlinearity $u\overline u$ in $H^{-1/4}(\R)$}

\author{Yuzhao Wang}
\address{LMAM, School of Mathematical Sciences, Peking University, Beijing
100871, China}

\email{wangyuzhao2008@gmail.com}

\begin{abstract}
In this note, we give another approach to the local well-posedness
of quadratic Schr\"odinger equation with nonlinearity $u\overline u$
in $H^{-1/4}$, which was already proved by Kishimoto \cite{kis}. Our
resolution space is $l^1$-analogue of $X^{s,b}$ space with low
frequency part in a weaker space $L^{\infty}_{t}L^2_x$. Such type
spaces were developed by Guo. \cite{G} to deal the KdV endpoint
$H^{-3/4}$ regularity.
\end{abstract}



\maketitle

\section{Introduction}

This paper is mainly concerned with the following equation
\begin{eqnarray}\label{eq:qsch}
\left\{
\begin{array}{l}
iu_t+u_{xx}=|u|^2,\quad u(x,t):\R\times \R\rightarrow \C,\\
u(x,0)=\phi(x)\in H^s(\R).
\end{array}
\right.
\end{eqnarray}
The low regularity for this equation was first studied by Kenig,
Ponce, Vega in \cite{KPV}, they proved the local well-posedness in
$H^{s}$, for $s>-1/4$, by using $X^{s,b}$ spaces. The local
well-posedness in $H^{-1/4}$ was already proved by Kishimoto
\cite{kis}, where Kishimoto solved \eqref{eq:qsch} in the spaces
$$
Z=X^{-1/4,1/2+\beta}+Y
$$
and
$$
Y=\{f\in \mathcal{S}'(\R^2);
\|f\|_{Y}=\|\jb{\xi}^{-1/4}\jb{\tau-\xi^2}^{3\beta}\hat{f}\|_{L^2_\xi
L^p_\tau}+\|\jb{\xi}^{1/4-2\beta}\jb{\tau-\xi^2}^{3\beta}\hat{f}\|_{L^2_\xi
L^2_\tau}\},
$$
with $0<\beta\leq 1/24$, $2\beta<1/p'<3\beta$ and $1/p+1/p' =1$.

We give another approach based on the argument developed by Guo.
\cite{G}, which solved the global well-posedness for KdV equation in
$H^{-3/4}$. Our resolution space is $l^1$-analogue of $X^{s,b}$
space with low frequency part in a weaker space
$L^{\infty}_{t}L^2_x$, so as a resolution space, it has simple form.

It is well known that $X^{s,b}$ failed for \eqref{eq:qsch} in
$H^{-1/4}$ because of the logarithmic divergences from $high\times
high \rightarrow low $ interactions, it is natural to use the weaker
structure in low frequency. We use $L^{\infty}_{t}L^{2}_{x}$ to
measure the low frequency part, however in \cite{G} Guo used
$L^{2}_{x}L^{\infty}_{t}$. The reason for this is that in the KdV
case, the $high\times low$ interactions has one derivative, and the
smoothing effect norm $L^\infty_xL^2_t$ was needed to absorb it.
This method can also be adapted to other similar problems where some
logarithmic divergences appear in the high-high interactions.

\begin{theorem}\label{thmlwp} The initial value problem
\eqref{eq:qsch} is local well-posedness in $H^{-1/4}$.
\end{theorem}

For $f\in \Sch'$ we denote by $\widehat{f}$ or $\ft (f)$ the Fourier
transform of $f$. We denote by $\ft_x$ the Fourier transform on
spatial variable. Let $\mathbb{Z}$ and $\mathbb{N}$ be the sets of
integers and natural numbers respectively, $\Z_+=\N \cup \{0\}$. For
$k\in \Z_+$ let ${I}_k=\{\xi: |\xi|\in [2^{k-1}, 2^{k+1}]\}, \ k\geq
1;$ $I_0=\{\xi: |\xi|\leq 2\}$. Let $\eta_0: \R\rightarrow [0, 1]$
denote an even smooth function supported in $[-8/5, 8/5]$ and equal
to $1$ in $[-5/4, 5/4]$. We define $\psi(t)=\eta_0(t)$. For $k\in
\Z$ let $\eta_k(\xi)=\eta_0(\xi/2^k)-\eta_0(\xi/2^{k-1})$ if $k\geq
1$ and $\eta_k(\xi)\equiv 0$ if $k\leq -1$. For $k\in \Z_+$, define
$P_k$
 by
$ \widehat{P_ku}(\xi)=\eta_k(\xi)\widehat{u}(\xi). $ For $l\in \Z$
let $ P_{\leq l}=\sum_{k\leq l}P_k, P_{\geq l}=\sum_{k\geq l}P_k. $

For $u_0\in \Sch'(\R)$, we denote $W(t)u_0=e^{it\partial_x^2}u_0$
 defined by $
\ft_x(W(t)\phi)(\xi)=\exp[-i\xi^2t]\widehat{\phi}(\xi) $.

For $k\in \Z_+$ we define the dyadic $X^{s,b}$-type normed spaces
$X_k=X_k(\R^2)$,
\begin{eqnarray}
X_k=\left\{f\in L^2(\R^2):
\begin{array}{l}
f(\xi,\tau) \mbox{ is supported in } I_k\times\R \mbox{ and }\\
\norm{f}_{X_k}=\sum_{j=0}^\infty
2^{j/2}\norm{\eta_j(\tau+\xi^2)\cdot f}_{L^2},
\end{array}
\right\}
\end{eqnarray}
thus we have $\|\widehat{f}\|_{L^{1}_{\tau}L^{2}_{\xi}}\leq
\|f\|_{X_k}$. For $-3/4\leq s\leq 0$, we define our resolution
spaces
\begin{eqnarray}\label{barF}
\bar{F}^s=\{u\in \Sch'(\R^2):\norm{u}_{\bar{F}^s}^2=\sum_{k \geq
1}2^{2sk}\norm{\eta_k(\xi)\ft(u)}_{X_k}^2+\norm{P_{\leq
0}(u)}_{_{L_t^\infty L_x^2}}^2<\infty\}.
\end{eqnarray}
It is easy to see that  for $k\in \Z_+$£¬
\begin{eqnarray}\label{em}
\norm{P_k(u)}_{L_t^\infty L_x^2}\les \norm{\ft[P_k(u)]}_{X_k},
\end{eqnarray}
as a consequence, we have $\norm{u}_{L_t^\infty H^s}\les
\norm{u}_{\bar{F}^s}$.

Let $a_1, a_2, a_3\in \R$, define $a_{max}=\max{\{a_1,a_2,a_2\}}$,
same as $a_{min}, a_{med}$. Usually we use $k_1,k_2,k_3$ and
$j_1,j_2,j_3$ to denote integers, $N_i=2^{k_i}$ and $L_i=2^{j_i}$
for $i=1,2,3$ to denote dyadic numbers.

\section{Dyadic Bilinear Estimates}

In this section we will give some dyadic bilinear estimates for next
section. We define
$$
D_{k,j}=\{(\xi,\tau): \xi \in [2^{k-1},2^{k+1}] \mbox{ and }
\tau+\xi^2\in I_j\}, \quad k\in \Z, j\in \Z_+.
$$
Following the $[k;Z]$ methods \cite{Taokz} the bilinear estimates in
$X^{s,b}$ space reduce to some dyadic summations: for any
$k_1,k_2,k_3\in \Z$ and $j_1,j_2,j_3\in \Z_+$
\begin{eqnarray}\label{eq:3zmult}
\sup_{(u_{k_2,j_2},\ v_{k_3,j_3})\in
E}\norm{1_{D_{k_1,j_1}}(\xi,\tau)\cdot
u_{k_2,j_2}*v_{k_3,j_3}(\xi,\tau)}_{L_{\xi,\tau}^2}
\end{eqnarray}
where $E=\{(u,v):\norm{u}_2,\ \norm{v}_2\leq 1 \mbox{ and } \supp(u)
\subset D_{k_2,j_2},\ \supp(v) \subset \widetilde{D}_{k_3,j_3}\}$
and $\widetilde{D}_{k_3,j_3}=\{(\xi,\tau); (-\xi, -\tau)\in
D_{k_3,j_3}\}$. By checking the support properties, we get that in
order for \eqref{eq:3zmult} to be nonzero one must have
\begin{align}
 |k_{max}-k_{med}|\leq 3 , \text{ and } j_{max}\geq k_{max}+k_{min}-10
\end{align}
The following sharp estimates on \eqref{eq:3zmult} were obtained
in \cite{Taokz}.

\begin{lemma}[Proposition 11.1, \cite{Taokz} (--++) case]\label{pchar}
Let $k_1,k_2,k_3 \in \Z$ and $j_1,j_2,j_3\in \Z_+$. Let
$N_i=2^{k_i}$ and $L_i=2^{j_i}$ for $i=1,2,3$. Then

(i) If $N_{max}\sim N_{min}$ and $L_{max}\sim N_{max}N_{min}$, then
we have
\begin{equation}\label{eq:chari}
\eqref{eq:3zmult} \les L_{min}^{1/2}L_{med}^{1/4}.
\end{equation}

(ii) If $N_1\sim N_3 \gg N_2$ and $N_{max}N_{min}\sim L_2=L_{max}$,
and $N_1\sim N_2 \gg N_3$ and $N_{max}N_{min}\sim L_3=L_{max}$, then
\begin{equation}\label{eq:charii}
\eqref{eq:3zmult} \les L_{min}^{1/2}L_{med}^{1/2}N_{min}^{-1/2}.
\end{equation}

(iii) In all other cases, we have
\begin{equation}\label{eq:chariii}
\eqref{eq:3zmult} \les
L_{min}^{1/2}N_{max}^{-1/2}\min(N_{max}N_{min},L_{med})^{1/2}.
\end{equation}
\end{lemma}

\section{Proof of Theorem \ref{thmlwp}}
For $u,v\in \bar{F}^s$ we define the bilinear operator
\begin{eqnarray}
B(u,v)=\psi\big(\frac{t}{4}\big)\int_0^tW(t-\tau)\partial_x\big(\psi^2(\tau)u(\tau)\cdot
v(\tau)\big)d\tau.
\end{eqnarray}
As in \cite{G}, the proof for Theorem \ref{thmlwp} reduce to showing
the boundness of $B:\bar{F}^{-1/4}\times \bar{F}^{-1/4}\rightarrow
\bar{F}^{-1/4}$.
\begin{lemma}[Linear estimates]\label{proplineares}
(a) Assume $s\in \R$,  $\phi \in H^s$. Then there exists $C>0$ such
that
\begin{eqnarray}
\norm{\psi(t)W(t)\phi}_{\bar{F}^s}\leq C\norm{\phi}_{H^{s}}.
\end{eqnarray}

(b) Assume $s\in \R, k\in \Z_+$ and $(i+\tau-\xi^3)^{-1}\ft(u)\in
X_k$. Then there exists $C>0$ such that
\begin{eqnarray}
\normo{\ft\left[\psi(t)\int_0^t W(t-s)(u(s))ds\right]}_{X_k}\leq
C\norm{(i+\tau-\xi^3)^{-1}\ft(u)}_{X_k}.
\end{eqnarray}
\end{lemma}

\begin{proof}
Such linear estimates have appeared in many literatures, see for
example \cite{In-Ke}.
\end{proof}

\begin{lemma}[Bilinear estimates]\label{propbilinearbd}
Assume $-1/4\leq s\leq 0$. Then there exists $C>0$ such that
\begin{eqnarray}\label{eq:bilinearbd}
\norm{B(u,v)}_{\bar{F}^s}\leq
C(\norm{u}_{\bar{F}^s}\norm{v}_{\bar{F}^{-1/4}}+\norm{u}_{\bar{F}^{-1/4}}\norm{v}_{\bar{F}^s})
\end{eqnarray}
hold for any $u,v\in \bar{F}^s$.
\end{lemma}
\begin{proof} It is easy to see
\begin{align*}
\norm{B(u,v)}_{\bar{F}^s}\les& \norm{P_{\geq 1}B(P_{\geq 1}u,P_{\geq
1}v)}_{F^s}+\norm{P_{\geq 1}B(P_{\geq 1}u,P_{
0}v)}_{\bar{F}^s}\\&+\norm{P_{\geq 1}B(P_{0}u,P_{\geq
1}v)}_{\bar{F}^s}+\norm{P_{\geq 1}B(P_{0}u,P_{
0}v)}_{\bar{F}^s}+\norm{P_{0}B(u,v)}_{\bar{F}^s} \\\triangleq&
A+B+C+D+E
\end{align*}
We notice that there is no low frequency in part $A$, so the proof
for part A do not involve the special structure in low frequency,
and standard $X^{s,b}$ argument will suffice, we omit the proof.

The proof for part B, C and D are similar, we just consider part B
for example. By definition and Lemma \ref{proplineares} (b), let
$S_B=\{(k_1,k_3); k_1,k_3 \geq 1, |k_{1}-k_{3}|\leq 5\}$, then
\begin{align}\label{b}
B^2\les& \sum_{(k_1,k_3)\in S_B} 2^{2sk_3}\Big(\sum_{j_3\geq
0}2^{-j_3/2}\norm{1_{D_{k_3,j_3}}\widehat{\psi(t)P_{k_1}u}*\widehat{P_{0}v}}_{L_{\xi,\tau}^2}\Big)^2\nonumber\\
\les& \sum_{(k_1,k_3)\in S_B}
2^{2sk_3}\|\psi(t)P_{k_1}u\|^2_{L^2}\|P_{0}v\|^2_{L^{\infty}}
\les\sum_{(k_1,k_3)\in S_B}
2^{2sk_3}\|P_{k_1}u\|^2_{L^{\infty}_{t}L_x^2}\|P_{0}v\|^2_{L^{\infty}}
\end{align}
which is sufficient by Bernstein inequality and \eqref{em}.

Now we turn to part D.  Denote $Q(u,v)=P_{\leq 0}B(P_{k_1}u,
P_{k_2}\bar{v})$. By straightforward computations,
\begin{eqnarray*}
\ft\left[Q(u,\bar{v})\right](\xi,\tau)=c\int_\R
\frac{\widehat{\psi}(\tau-\tau')-\widehat{\psi}(\tau+\xi^2)}{\tau'+\xi^2}\eta_0(\xi)\int_Z
\wh{P_{k_1}u}(\xi_1,\tau_1)\wh{P_{k_2}\bar{v}}(\xi_2,\tau_2)\
d\tau'.
\end{eqnarray*}
where $Z=\{\xi=\xi_1+\xi_2, \tau'=\tau_1+\tau_2\}$. Fixing $\xi \in
\R$, we decomposing the hyperplane as following
\begin{eqnarray*}
A_1&=&\{\xi=\xi_1+\xi_2,\tau'=\tau_1+\tau_2: |\xi|\les 2^{-k_1}\};\\
A_2&=&\{\xi=\xi_1+\xi_2,\tau'=\tau_1+\tau_2: |\xi|\gg 2^{-k_1},\\&&
|\tau_1+\xi_1^2|\ll 2^{k_1}|\xi|, |\tau_2-\xi_2^2|\ll  2^{k_1}|\xi|\};\\
A_3&=&\{\xi=\xi_1+\xi_2,\tau'=\tau_1+\tau_2: |\xi|\gg 2^{-k_1},
|\tau_1+\xi_1^2|\ges  2^{k_1}|\xi|\};\\
A_4&=&\{\xi=\xi_1+\xi_2,\tau'=\tau_1+\tau_2: |\xi|\gg 2^{-k_1},
|\tau_2-\xi_2^2|\ges  2^{k_1}|\xi|\}.
\end{eqnarray*}
Then we get
\[\ft\left[Q(u,\bar{v})\right](\xi,\tau)=I+II+III,\]
where
\begin{align*}
I=&C\int_\R
\frac{\widehat{\psi}(\tau-\tau')-\widehat{\psi}(\tau+\xi^2)}{\tau'+\xi^2}\eta_0(\xi)
\int_{A_1}\wh{P_{k_1}u}(\xi_1,\tau_1)\wh{P_{k_2}\bar{v}}(\xi_2,\tau_2)d\tau',\\
II=&C\int_\R
\frac{\widehat{\psi}(\tau-\tau')-\widehat{\psi}(\tau+\xi^2)}{\tau'+\xi^2}\eta_0(\xi)
\int_{A_2}\wh{P_{k_1}u}(\xi_1,\tau_1)\wh{P_{k_2}\bar{v}}(\xi_2,\tau_2)d\tau',\\
III=&C\int_\R
\frac{\widehat{\psi}(\tau-\tau')-\widehat{\psi}(\tau+\xi^2)}{\tau'+\xi^2}\eta_0(\xi)
\int_{A_3\cup
A_4}\wh{P_{k_1}u}(\xi_1,\tau_1)\wh{P_{k_2}\bar{v}}(\xi_2,\tau_2)d\tau'.
\end{align*}

We consider first the the term $I$. By \eqref{em} and Proposition
\ref{proplineares} (b),
\[\norm{\ft^{-1}(I)}_{L_t^\infty L_x^2}\les \norm{I}_{X_0}\les \Big\|(i+\tau'+\xi^2)^{-1}\eta_0(\xi)
\int_{A_1}\wh{P_{k_1}u}(\xi_1,\tau_1)\wh{P_{k_2}\bar{v}}(\xi_2,\tau_2)\Big\|_{X_0},\]
since in the set $A_1$ we have $|\xi|\les 2^{-k_1}$, thus we
continue with
\begin{eqnarray*}
\les \sum_{k_3\leq -k_1+10}\sum_{j_3\geq 0} 2^{-j_3/2}\sum_{j_1\geq
0, j_2\geq 0}\norm{1_{D_{k_3,j_3}}\cdot
u_{k_1,j_1}*v_{k_2,j_2}}_{L^2}
\end{eqnarray*}
where
\begin{eqnarray}
u_{k_1,j_1}=\eta_{k_1}(\xi)\eta_{j_1}(\tau+\xi^2)\widehat{u},\
v_{k,j_2}=\eta_k(\xi)\eta_{j_2}(\tau-\xi^2)\widehat{\overline{v}}.
\end{eqnarray}
Using Proposition \ref{pchar} (iii), then we get
\begin{eqnarray*}
\norm{\ft^{-1}(I)}_{L_t^\infty L_x^2}&\les& \sum_{k_3\leq
-k_1+10}\sum_{j_i\geq 0} 2^{-j_3/2}2^{j_{min}/2}2^{k_3/2}
\norm{u_{k_1,j_1}}_{L^2}\norm{v_{k_2,j_2}}_{L^2}\\
&\les&
2^{-k_1/2}\norm{\wh{P_{k_1}u}}_{X_{k_1}}\norm{\wh{P_{k_2}v}}_{X_{k_2}},
\end{eqnarray*}
which suffices to give the bound for the term $I$ since
$|k_1-k_2|\leq 5$.

Next we consider the contribution of the term $III$. As term $I$, by
\eqref{em} and Proposition \ref{proplineares} (b),
\begin{eqnarray*}
\norm{\ft^{-1}(III)}_{L_t^\infty L_x^2}&\les&
\normo{(i+\tau'+\xi^2)^{-1}\eta_0(\xi)\int_{A_3\cup
A_4}\wh{P_{k_1}u}(\xi_1,\tau_1)\wh{P_{k_2}\overline{v}}(\xi_2,\tau_2)}_{X_0}\\
&\les&\sum_{-k_1\leq k_3\leq 0}\sum_{j_3\geq 0}
2^{-j_3/2}\sum_{j_1\geq 0, j_2\geq 0}\norm{1_{D_{k_3,j_3}}\cdot
u_{k_1,j_1}*v_{k_2,j_2}}_{L^2}.
\end{eqnarray*}
Without loss of generality, we assume $|\tau_1+\xi_1^2|\ges
|\xi\xi_1|$, applying Proposition \ref{pchar} (iii), then we get
\begin{eqnarray*}
\norm{\ft^{-1}(III)}_{L_t^\infty L_x^2}&\les& \sum_{-k_1\leq k_3\leq
0}\sum_{j_1\geq k_3+k_1-10, j_2\geq 0} 2^{j_{2}/2}2^{-k_1/2}
\norm{u_{k_1,j_1}}_{L^2}\norm{v_{k_2,j_2}}_{L^2}\\
&\les&
2^{-k_1/2}\norm{\wh{P_{k_1}u}}_{X_{k_1}}\norm{\wh{P_{k_2}u}}_{X_{k_2}},
\end{eqnarray*}
which suffices to give the bound for the term $III$ since
$|k_1-k_2|\leq 5$.

Now we consider the main contribution term: term $II$.
 By direct computation, we get
\begin{align*}
\ft_t^{-1}(II)=\psi(t)\int_0^t e^{-i(t-s)\xi^2}\eta_0(\xi)i\xi
\int_{\R^2}e^{is(\tau_1+\tau_2)} \
\int_{\xi=\xi_1+\xi_2}{u_{k_1}}(\xi_1,\tau_1){v_{k_2}}(\xi_2,\tau_2)\
d\tau_1 d\tau_2ds
\end{align*}
where
\begin{align*}
u_{k_1}(\xi_1,\tau_1)=\eta_{k_1}(\xi_1)1_{\{|\tau_1+\xi_1^2|\ll
 2^{k_1}|\xi|\}}\wh{u}(\xi_1,\tau_1), \text{ }
v_{k_2}(\xi_2,\tau_2)=\eta_{k_2}(\xi_2)1_{\{|\tau_2-\xi_2^2|\ll
2^{k_1}|\xi|\}}\wh{\overline{v}}(\xi_2,\tau_2).
\end{align*}
By a change of variable $\tau_1'=\tau_1+\xi_1^2$,
$\tau_2'=\tau_2-\xi_2^2$, we get
\begin{eqnarray*}
\ft_t^{-1}(II)&=&\psi(t)e^{-it\xi^2}\eta_0(\xi)\int_0^t e^{is\xi^2}
\int_{\R^2}e^{is(\tau_1+\tau_2)}\\
&&\times \
\int_{\xi=\xi_1+\xi_2}e^{-is\xi_1^2}{u_{k_1}}(\xi_1,\tau_1-\xi_1^2)e^{is\xi_2^2}{v_{k_2}}(\xi_2,\tau_2+\xi_2^2)\
d\tau_1 d\tau_2ds\\
&=&\psi(t)e^{-it\xi^2}\eta_0(\xi)\int_{\R^2}e^{it(\tau_1+\tau_2)}\int_{\xi=\xi_1+\xi_2}
\frac{e^{it(-\xi_1^2+\xi_2^2+\xi^2)}-e^{-it(\tau_1+\tau_2)}}{\tau_1+\tau_2-\xi_1^2+\xi_2^2+\xi^2}\\
&&\times \
{u_{k_1}}(\xi_1,\tau_1-\xi_1^2){v_{k_2}}(\xi_2,\tau_2+\xi_2^2)\
d\tau_1 d\tau_2.
\end{eqnarray*}
Then by Plancherel Theorem and H\"older inequality, we can bound $
\|\ft_t^{-1}(II)\|_{L_\xi}$ by
\begin{align*}
&\int_{\R^2}\Big\|\int_{\xi=\xi_1+\xi_2}\frac{\eta_0(\xi)}{|\tau_1+\tau_2-\xi_1^2+\xi_2^2+\xi^2|}
|
u_{k_1}(\xi_1,\tau_1-\xi_1^2)v_{k_2}(\xi_2,\tau_2+\xi_2^2)|\Big\|_{L^2_\xi}\
d\tau_1 d\tau_2\\
\les&\int_{\R^2}\sum_{-k_1\leq k\leq
0}2^{k/2}\Big\|\int_{\xi=\xi_1+\xi_2}\frac{\chi_k(\xi)}{|\xi \xi_1|}
|u_{k_1}(\xi_1,\tau_1-\xi_1^2)v_{k_2}(\xi_2,\tau_2+\xi_2^2)|\Big\|_{L^\infty_\xi}\
d\tau_1 d\tau_2\\
\les&
2^{-k_1/2}\|u_{k_1}\|_{L^{1}_{\tau_2}L_{\xi_2}^2}\|v_{k_2}\|_{L^{1}_{\tau_3}L_{\xi_3}^2}\les
2^{-k_1/2}\norm{\wh{P_{k_1}u}}_{X_{k_1}}\norm{\wh{P_{k_2}u}}_{X_{k_2}}.
\end{align*}
where we use $|\tau_1+\tau_2-\xi_1^2+\xi_2^2+\xi^2|\gtrsim |\xi
\xi_1|$, which completes the proof of the lemma.
\end{proof}

\noindent{\bf Acknowledgment.} The author is very grateful to
Professor Zihua Guo for encouraging the author to work on this
problem and helpful conversations.

\end{document}